\documentclass[12pt,intlimits]{amsart}
\usepackage{amssymb,amsmath,amsthm,upref,cite}

\numberwithin{equation}{section}

\theoremstyle{plain}
\newtheorem{theorem}{Theorem}[section]

\newtheorem{lemma}{Lemma}[section]

\theoremstyle{remark}
\newtheorem{remark}[theorem]{Remark}

\newcommand{\const}{\operatorname{const}}
\newcommand{\p}{\partial}

\begin{document}

\title{Lemniscates do not survive Laplacian growth}

\author[Khavinson]{D. Khavinson}
\address{Department of Mathematics \& Statistics, University of South Florida, Tampa, FL 33620-5700, USA}
\email{dkhavins@cas.usf.edu}

\author[Mineev-Weinstein]{M. Mineev-Weinstein}
\address{Los Alamos National Laboratory, MS-365, Los Alamos, NM 87545, USA}
%\curraddr{}
\email{mariner@lanl.gov}

\author[Putinar]{M. Putinar}
\address{Department of Mathematics, University of California, Santa Barbara, CA 93106, USA}
\email{mputinar@math.ucsb.edu}

\author[Teodorescu]{R. Teodorescu}
\address{Department of Mathematics \& Statistics, University of South Florida, Tampa, FL 33620-5700, USA}
\email{razvan@cas.usf.edu}

\thanks{Part of this work was done during the first and third
authors' visit to LANL and was supported by the LDRD project 20070483 ``Minimal Description of
Complex 2D Shapes" at LANL. Also, D. Khavinson and M. Putinar gratefully
acknowledge partial support by the National Science Foundation.}

\date{}

%\begin{abstract}
%\end{abstract}

\maketitle

\section{Introduction}\label{sec1}

Many moving boundary processes in the plane, e.g., solidification, electrodeposition, viscous fingering, bacterial growth, etc., can be mathematically modeled by the so-called Laplacian growth \cite{MPT, VE}. In a nutshell, it can be described by the equation
\begin{equation}
\label{eq1.1}
V(z)=\partial_n g_{\Omega(t)}(z,\zeta),
\end{equation}
where $V$ is the normal component of the velocity of the boundary
$\partial\Omega(t)$ of the moving domain
$\Omega(t)\subset\mathbb{R}^2\simeq\mathbb{C}$,
$z\in\partial\Omega(t)$, $t$ is time, $\dfrac\partial{\partial n}$
denotes the normal derivative on $\partial\Omega(t)$ and
$g_{\Omega(t)}(z,\zeta)$ is the Green function for the Laplace operator
in the domain $\Omega(t)$ with a unit source at the point
$\zeta\in\Omega(t)$. Equation \eqref{eq1.1} can be elegantly rewritten
as the area-preserving diffeomorphism
\begin{equation}
\label{eq1.2}
\Im\left(\bar{z}_tz_\theta\right)=1,
\end{equation}
where $\Im$ denotes the imaginary part of a complex number,  $\partial\Omega(t) := $ $\{z := z\left(t,\theta\right)\}$  is the
moving boundary parametrized by $w=e^{i\theta}$ on the unit circle
and the conformal mapping from, say, the exterior of the unit disk
$\mathbb{D}^+:=\{|w|>1\}$ onto $\Omega(t)$ with the normalization
$z(\infty)=\zeta, z'(\infty) > 0$.

The equation \eqref{eq1.2},  named {\emph{ Laplacian growth}} or the {\emph {Polubarinova - Galin}}
equation in modern literature, was first derived by Polubarinova-Kochina \cite{PK}
and Galin \cite{Galin} in 1945, as a description of secondary oil recovery processes.

This equation is known to be integrable \cite{MWZ}, and as such possesses an infinte number of conserved quantities. 
More precisely, it admits conserved moments $c_n=\int\limits_{\Omega(t)}z^ndx\,dy$, where $n$ runs over either all non-negative or all non-positive integers depending on whether domains $\Omega(t)$ are finite or infinite. At the same time \eqref{eq1.2} admits an impressive number of closed-form solutions.

For the background, history, generalizations, references, connections to the theory of quadrature domains and other branches of mathematical physics we refer the reader to \cite{EGKP,KMP,Sh,VE, MPT, MWZ} and the references therein.

In section \S \ref{sec2} of this paper, we show that any continuous
chain of polynomial lemniscates of order $n$:
$\Gamma_t:=\{\left|P(z,t)\right|=1\}$,
$P(z,t)=a(t)\prod\limits_{j=1}^n\left(z-\lambda_j(t)\right)$, where $a(t)$ is real-valued, is
destroyed instantly under the Laplacian growth process described in
\eqref{eq1.1}, with $\Omega(t)=\{|P(z,t)|>1\}$,
$\zeta=\infty$, unless $n=1$, $\lambda_1(t)=\const$ and
$\left\{\Gamma_t=\partial\Omega(t)\right\}$ is simply a family of
concentric circles. Here the roots $\lambda_j(t)$ of $P(z,t)$ are
all assumed to be inside $\Omega'_t:=\{|P(z,t)|<1\}$, so $\Omega$ and $\Omega'$
are simply connected.

This result shows that unlike quadrature domains (cf.
\cite{EGKP,Sh}) that are preserved under the Laplacian growth
process, lemniscates for which all  the roots of the defining polynomial are in $\Omega'_t$ are instantly destroyed,
 except for the trivial case of concentric circles. This, incidentally,  agrees with a
well-known fact --- cf. \cite{EKS} --- that lemniscates which are
also quadrature domains must be circles.
The proof of the theorem for the case of Laplacian growth is given in \S\ref{sec2}. 

In \S\ref{sec3} we extend the result of \S\ref{sec2} to all the growth processes 
that are invariant under time-reversal
and for which the boundary velocity is given by
\begin{equation} \label{def_process}
V (z)= \chi(z) \p_n g_{\Omega(t)}(z,\zeta),
\end{equation}
with $\chi(z)$ is a bounded, real, positive function on $\Gamma_t$.  
%%and subject to the condition
%\begin{equation}\label{boundary_normal}
%\vec { \nabla} \cdot \vec{V}(z, t)  = 0, \quad z \in U_t.
%\end{equation}
Invariance under time-reversal is defined here in the following way: if the boundary $\Gamma_{t+dt}$ is the image of $\Gamma_t$ under a map $f_{(t, dt)}:  z_{t}\in \Gamma_t \mapsto z_{t+dt} \in \Gamma_{t +dt}$, then $f_{(t+dt, -dt)}\circ f_{(t, dt)} = \mathbb{Id}$. 

We conclude with a few remarks in \S\ref{sec4}. 

\section{Destruction of Lemniscates}\label{sec2}

\begin{theorem}\label{thm2.1}
Suppose that a family of moving boundaries $\Gamma_{t}$,  (where
$t>0$ is time), produced by a Laplacian growth process, is  a family of
polynomial lemniscates $\{|P(z,t)|=1\}$, where
$P(z,t)=a(t)\prod\limits_{j=1}^n\left[z-\lambda_j(t)\right]$, 
and all $\lambda_j(t)$ are assumed to be inside $\Gamma_t$. Then,
$n=1$ and $\lambda_1=\const$, i.e., $\Gamma_t$ is a family of
concentric circles.
\end{theorem}

\begin{proof}
Let $\Omega_t=\{z:|P(z,t)|>1\}$, $\mathbb{D}^+=\{|w|>1\}$. The
function $\varphi(t): \Omega_t\to\mathbb{D}^+$,
$w=\varphi(z, t)=\sqrt[n]{P(z,t)}$, where we choose the branch for
the $n-$th root so that $\varphi'(t, \infty)>0$, maps $\Omega_t$
conformally onto $\mathbb{D}_+$, $\varphi(t, \infty)=\infty$. It is useful to note 
that on $\Gamma_t$,  $P(z,t)=w^n$, $|w|=1$ and does not depend on $t$. 
This is because for any two moments of time $t, \tau$, we have $w(t)(.) = w(\tau)\circ \kappa(t, \tau)(.)$, 
where $\kappa$ is a M\"{o}bius automorphism of the disk. In our case, 
$\kappa(t, \tau)(\infty) = \infty$, so $\kappa(t, \tau)(z) = e^{i \alpha} z, \alpha\in \mathbb{R}$, but
since it also fixes the argument at $\infty$, $\kappa$ is the identity. 
 
Therefore, we have (where, as is customary, we
denote the partial $t$-derivative by a ``dot''):
\begin{equation}
\label{eq2.1}
\dot{P}+P'_z\dot{z}=0.
\end{equation}
Since $\varphi(t)$ maps $\Gamma_t$ onto the unit circle, we have
$z(t)=\Psi(t, w)$, where $\Psi(t, w)=\varphi^{-1}(t, z)$. 
 We also have on $\Gamma_t$, by differentiating $P(z(w),t)=w^n$ with respect to $w$,
\begin{equation}
\label{eq2.2}
P'_z \cdot z_w = nw^{n-1}
\end{equation}
or
\begin{equation}
\label{eq2.3}
wz_w=\frac{nw^n}{P'_z}=\frac{nP}{P'_z}.
\end{equation}

From \eqref{eq2.1}, conjugating, we infer
\begin{equation}
\label{eq2.4}
\overline{\dot{z}}=-\frac{\overline{\dot{P}}}{\overline{P'_z}}.
\end{equation}
Parametrize the unit circle by $w=e^{i\theta}$, $0\le\theta\le
2\pi$. Then, from \eqref{eq2.3}, it follows that
  we have on $\Gamma_t$ (since $(z(t)=z(w,t)=z(w(\theta,t)))$),
\begin{equation}
\frac 1i\,z_\theta:=\frac{\partial z}{i\partial\theta}=z_ww=\frac{nP}{P'_z}.\label{eq2.5}
\end{equation}
Combining \eqref{eq2.4} and \eqref{eq2.5} yields ($\Re$ stands for the real part):
\begin{equation}
\label{eq2.6}
\Re\left(\overline{\dot{z}}\frac1i\,z_\theta\right)
=\Re\left(-\frac{\overline{\dot{P}}}{\overline{P'_z}}
\cdot\frac{nP}{P'_z}\right).
\end{equation}
Also,
\begin{equation}
\label{eq2.7} \Re\left(\overline{\dot{z}}\frac{\partial
z}{i\partial\theta}\right)
=\Im\left(\overline{\dot{z}}z_{\theta}\right)=1,
\end{equation}
where in the last equality we used the hypothesis that the
lemniscates $\Gamma_t:=\{z(t,\theta)\}$ satisfy the main equation
\eqref{eq1.2} of Laplacian growth processes--- cf. \cite[\S4]{KMP}. Hence, \eqref{eq2.7}, \eqref{eq2.4} and
\eqref{eq2.5} imply that
\begin{equation}
\label{eq2.8}
\frac1n\Re\left(\frac{\dot{P}}{P'_z}\,\overline{iz_\theta}\right)
=\Re\left(\frac{\dot{P}}{P'_z}\;
\frac{\overline{P}}{\overline{P'_z}}\right)=\frac{-1}n.
\end{equation}
Or, we can rewrite \eqref{eq2.8} as
\begin{equation}
\label{eq2.9}
\Re\left(\dot{P}\overline{P}\right)
=-\frac1n\left|P'_z\right|^2.
\end{equation}
Thus, we are finally arriving at
\begin{equation}
\label{eq2.10} \frac
d{dt}\left(|P|^2\right)=-\frac{1}{2n}\left|P'_z\right|^2.
\end{equation}
Therefore,  \eqref{eq2.10} holds on the
 lemniscates $\Gamma_t=\{|P(z,t)|=1\}$ that are assumed to be interfaces of a Laplacian growth process. Now the theorem follows from the
following.

\begin{lemma}
Let $t$ be the time variable, $P(z,t)=a(t)\prod\limits_1^n\left(z-\lambda_i(t)\right),$ be a ``flow'' of $n$-degree polynomials. 
Assume that the lemniscates $\Gamma_t:=\{|P(z,t)|=1\}$ all have connected interiors $\{|P(z,t)|<1\}$ and a generalized 
equation \eqref{eq2.10} holds on $\Gamma_t$; i.e.,
\begin{equation}
\label{eq2.11} \frac d{dt}\left(|P(z,t)|^2\right)
-c(t)\left|P'_z(z,t)\right|^2=0,
\end{equation}
where the function $c(t)$ is real-valued, depends on $t$ only and, hence,  is a
constant on $\Gamma_t$. Then, $n=1$, $\lambda_1=\lambda_1(t)=\const$
and $\Gamma_t$ is a family of concentric circles centered at
$\lambda_1$.
\end{lemma}

\noindent \emph{Proof of the Lemma.} 
Our hypothesis implies that all  polynomials $|P(z,t)|^2-1$ are irreducible. Hence, using Hilbert's Nullstellensatz (e.g., c.f. \cite{Bourbaki}, Proposition 3.3.2), we infer from \eqref{eq2.11} that
\begin{equation}
\label{eq2.12} \frac d{dt}\left(|P(z,t)|^2\right)
-c(t)\left|P'_z(z,t)\right|^2 =B(t)\left(|P(z,t)|^2 -1\right).
\end{equation}
Equation \eqref{eq2.12} holds for all $z\in\mathbb{C}$ and for an interval of time $t$, and for each $t$, both sides are real-analytic functions in $z$ and $\overline{z}$. Hence, we can ``polarize'' \eqref{eq2.12}, i.e., replace $\overline{z}$ by an independent complex variable $\xi$. (This is due to a simple observation: real-analytic functions of two variables are nothing else but restrictions of holomorphic functions in $z,\xi$-variables to the plane $\left\{\xi=\overline{z}\right\}$. Hence, if two real-analytic functions coincide on that plane, they coincide in $\mathbb{C}^2$ as well.) Denoting by $P^{\#}$ the polynomial whose coefficients are obtained from $P$ by complex conjugation, we have \eqref{eq2.12} in a ``polarized'' form holding for $(z,\xi)\in\mathbb{C}^2$:
\begin{multline}
\label{eq2.13}
\frac d{dt}\left(P(z,t)P^{\#}(\xi,t)\right)-c(t)\left(P'_z(z,t)\cdot\left(P^{\#}\right)'_{\xi}(\xi,t)\right) \\
=B(t)\left(P(z,t)P^{\#}(\xi,t)-1\right).
\end{multline}
Now let us denote by $k_j$ the multiplicity of the root $\lambda_j(t)$ of the polynomial $P(z, t)$, so that there are $m\le n$ distinct roots and $\sum_{j=1}^m k_j = n$.  Since
\begin{gather*}
P(z,t)=a(t)\prod_1^m\left(z-\lambda_j(t)\right)^{k_j}, \\
P^{\#}(\xi,t)=\bar{a}(t)\prod_1^m\left(\xi-\overline{\lambda_j(t)}\right)^{k_j},
\end{gather*}
dividing by $P(z,t)P^{\#}(\xi,t)$ we obtain:
\begin{multline}
\label{eq2.14} 2{\Re}\left ( \frac{\dot{a}}{a}\right )-\sum_1^m
\left(\frac{k_j \dot{\lambda}_j(t)}{z-\lambda_j(t)}
+\frac{k_j \overline{ \dot{\lambda}_j}(t)}{\xi-\overline{\lambda_j(t)}}\right) \\
-c(t)\left[\sum_1^m\frac{k_j}{z-\lambda_j(t)}\right] \cdot
\left[\sum_1^m\frac{k_j}{\xi-\overline{\lambda_j(t)}}\right] \\
=B(t)\left ( 1 - \frac{1}{P(z,t)P^{\#}(\xi,t)}\right ).
\end{multline}
Integrating \eqref{eq2.14} along a small circle centered at $\lambda_j(t)$, so that it does not enclose other zeros of $P$, yields for all $\xi$:
\begin{equation}
\label{eq2.15}
-k_j \dot{\lambda}_j(t)-c(t)\left(\sum_1^m\frac{k_i k_j}{\xi-\bar{\lambda}_i(t)}\right)
=-\frac{B(t)}{P^{\#}(\xi,t)}q_j,
\end{equation}
where $q_j = \frac{1}{(k_j-1)!}\left ( \frac{\p }{\p z}\right )^{k_j-1} \left [\frac{(z-\lambda_j)^{k_j}}{P(z,t)}\right ]_{z=\lambda_j}$. Letting $\xi\to\infty$ in \eqref{eq2.15} implies that $\dot{\lambda}_j(t)=0$ for all $j=1,\dotsc,n$. In other words, the ``nodes'' $\lambda_j(t)$ of all the lemniscates $\Gamma_t$ are fixed, i.e. do not move with time. So,
\begin{equation}
\label{eq2.16}
P(z,t)=a(t)\prod_1^n\left(z-\lambda_j\right)=a(t)\,Q(z).
\end{equation}
Substituting \eqref{eq2.16} into \eqref{eq2.13}, we obtain
\begin{multline}
\label{eq2.17}
\frac d{dt}\left(|a|^2\right)Q(z)Q^{\#}(\xi)-c(t)|a|^2Q'_z\left(Q^{\#}\right)'_\xi \\
=B(t)\left(|a|^2Q(z)Q^{\#}(\xi)-1\right).
\end{multline}
Comparing the leading terms (i.e., the coefficients at $z^n\xi^n$) in \eqref{eq2.17} yields
\begin{equation}
\label{eq2.18}
\frac d{dt}\left(|a|^2\right)=B(t)\,|a|^2.
\end{equation}
Therefore,
\begin{equation}
\label{eq2.19} c(t)\,|a|^2Q'_z\left(Q^{\#}\right)'_{\xi}=B(t),
\end{equation}
and thus $\deg Q'_z=0$, i.e., $n=\deg P=1$. The proofs of the Lemma and the Theorem are now complete.
\end{proof}

\section{Extending the theorem to growth processes invariant under time reversal} \label{sec3}

First, let us note that  any boundary $\Gamma_{t}$ is an equipotential line of the logarithmic potential 
\begin{equation} \label{potential}
\Phi(z) = \log |P_n(z, \lambda_i(t))|^2. 
\end{equation}
The boundary velocity of   the   general growth process defined in \eqref{def_process} can now be expressed as $\vec{V}(z) = \chi(z) \vec{\nabla} \Phi, z\in \Gamma_t, \chi(z) \in \mathbb{R}_+$.

As indicated in the Introduction, invariance under time-reversal is defined here in the following way: if the boundary $\Gamma_{t+dt}$ is the image of $\Gamma_t$ under a map $f_{(t, dt)}:  z_{t}\in \Gamma_t \mapsto z_{t+dt} \in \Gamma_{t +dt}$, then $f_{(t+dt, -dt)}\circ f_{(t, dt)} = \mathbb{Id}$. That means that the normal at $z_{t+dt} \in \Gamma_{t+dt}$ must be parallel to the normal at $z_t \in \Gamma_t$, which shows that $\Gamma_{t+dt}$ is perpendicular at every point to gradient lines of $\Phi$, and is therefore a level line of $\Phi$.  The displacement of the point $z_t$ becomes 
$$
z_{t+dt} -z_t = \chi(z) \vec{\nabla}\Phi(z_t) dt.
$$ 
Denoting by 
$$
\vec{E} = \vec{\nabla} \Phi = 2 \bar{\p} \Phi = 2\cdot \overline{\frac{P'_n(z, \lambda_i(t))}{P_n(z, \lambda_i(t))}}
$$ 
the gradient of the logarithmic potential and by $\vec{r} = z_t$, conservation of the normal (or gradient) direction becomes 
$$
\vec{E} (\vec{r} + \chi \vec{E}(\vec{r}) dt ) = \mu(z) \vec{E} (\vec{r}), 
$$
where $ \mu(z) = 1 +  m(z) dt,  m(z)=O(1),  m(z)  \in \mathbb{R}$, so after expanding in the infinitesimal time interval $dt$, 
%$$
%\vec{E}(\vec{r}) + \chi(z) dt (\vec{E} \cdot \vec{\nabla}) \vec{E}(\vec{r}) =\mu(z) \vec{E}(\vec{r}), 
%$$
%Equivalently,
$$
(\vec{E} \cdot \vec{\nabla}) \vec{E}(\vec{r})  = \frac{m(z) }{\chi(z)} \vec{E}(\vec{r}).
$$

\begin{remark}
The proportionality relation indicated above carries also the following physical significance: the dynamical system that we study is of \emph{frictional} type, where the \emph{acceleration} field (proportional to the force, or gradient of Green's function) is also proportional to the \emph{velocity}. In other words, the transport derivative (or Lie derivative) of the velocity field must be parallel to the velocity itself:
$$
\mathcal{L}_{\vec{V}} \vec{V} =[ {i}_{\vec{V}} \circ \textrm{d} - \textrm{d} \circ {i}_{\vec{V}} ]\vec{V}  = (\vec{V}\cdot \vec{\nabla})\vec{V} = \chi [\chi(\vec{E}\cdot \vec{\nabla})\vec{E} 
+  (\vec{E}\cdot \vec{\nabla \chi})\vec{E} ] 
$$
is parallel to $\vec{V}$ and therefore, to $\vec{E}$. 
\end{remark}

In complex notation, using the fact that $ (\vec{E} \cdot \vec{\nabla})  = \bar{E} \bar \p + E \p$, we obtain
$$
\overline{\frac{P'_n(z, \lambda_i(t))}{P_n(z, \lambda_i(t))}} = \delta(z) \frac{P'_n(z, \lambda_i(t))}{P_n(z, \lambda_i(t))} \cdot
\overline{ \left( \frac{P'_n(z, \lambda_i(t))}{P_n(z, \lambda_i(t))}\right )'}, \quad \delta(z) \in \mathbb{R},
$$
which (after multiplying both sides by $E(z)$) reduces to 
$$
\left [ {\frac{P'_n(z, \lambda_i(t))}{P_n(z, \lambda_i(t))}}\right ] ^{-2} { \left( \frac{P'_n(z, \lambda_i(t))}{P_n(z, \lambda_i(t))}\right )' } \in \mathbb{R},
$$
or 
\begin{equation} \label{cond}
\Im \left \{ \left [ \frac{P_n(z, \lambda_i(t))}{P'_n(z, \lambda_i(t))} \right ] ' \right \} = 0 , \quad (\forall) z \in \Gamma_{t}.
\end{equation}
We note that, since $E(z) = 2\bar{P'_n}/\bar P_n$ is the gradient of the Green's function for $\Omega_t$ and $\Omega_t$ is simply connected, 
it cannot vanish anywhere in $\Omega_t \cup \Gamma_t$, so all the  zeros of $P'_n(z)$, denoted by $\xi_k, k = 1, 2, \ldots, n-1$, are found inside the domain $\Omega'_t$. Then 
$$
\left [ \frac{P_n(z, \lambda_i(t))}{P'_n(z, \lambda_i(t))} \right]' = \frac{1}{n} + \sum_{k=1}^{n-1} \frac{A_k}{(z-\xi_k)^2},  \quad \xi_k \in \Omega'_t, 
$$
with $A_k$ constants. The imaginary part of this expression coincides with the imaginary part of an analytic function in $\Omega_t$, that is bounded there, so the condition \eqref{cond} can only be satisfied if the function is a constant. Since at $z \to \infty$ it vanishes, it follows that 
$$
\left [ \frac{P_n(z, \lambda_i(t))}{P'_n(z, \lambda_i(t))} \right ] '  = \frac{1}{n},
$$
which means that boundaries $\Gamma_t$ can only be concentric circles. 

\section{Concluding Remarks}\label{sec4}

\begin{enumerate}
\item  It is plausible that the result can be extended to rational lemniscates $\Gamma_t:=\left\{\left|R\left(z, t\right)\right|=1\right\}$,
where $R\left(z, t\right)$ are rational functions of degree $n$
where all the zeros are inside $\Gamma$, while all poles are in the
unbounded component of $\mathbb{C}\setminus\Gamma_t$.
\item  It is well-known that arbitrary ``shapes", i.e. Jordan curves can be arbitrarily close approximated by both lemniscates (Hilbert's theorem -- cf. \cite{Walsh}) and quadrature domains \cite{Bell}. At the same time our results imply that there are fundamental differences between these two classes of curves. We think it is interesting to pursue these observations in greater depth. 
\item From the argument in \S \ref{sec3} we can extract more. Suppose a family of Jordan curves, $\{ \Gamma_t \}_{t > 0}$ evolves by the flow along the velocity field $V(z)$ according to \eqref{def_process}. Assuming 
the invariance under time-reversal, the argument of \S\ref{sec3} can be used to prove that $\chi = $ const, i.e. the process is that of Laplacian growth. Invoking now well-known results on standard Hele-Shaw flows, we can at once conclude, e.g., that the process \eqref{def_process} continues for all times $t > 0$, i.e., the curves $\{ \Gamma_t \}$ move out to infinity such that $\cup_{t>t_0} \Gamma_t = {\mathbb{C}} \setminus \overline{\Omega}_{t_0}$, if and only if the initial curve $\Gamma_0$ is an ellipse and all the curves $\{ \Gamma_t \}$ are also ellipses homotetic with $\Gamma_0$ - cf. \cite{DiF}, also cf. \cite{FS}. 
\end{enumerate}

\providecommand{\bysame}{\leavevmode\hbox to3em{\hrulefill}\thinspace}
\providecommand{\MR}{\relax\ifhmode\unskip\space\fi MR }
\providecommand{\MRhref}[2]{%
  \href{http://www.ams.org/mathscinet-getitem?mr=#1}{#2}
}
\providecommand{\href}[2]{#2}

%\bibliographystyle{amsplain}
%\bibliography{lemniscates}

\end{document}